\documentclass[preprint,showkeys,preprintnumbers,amsmath,amssymb]{revtex4}
\usepackage{amsmath}
\usepackage{graphicx}
\usepackage{dcolumn}
\usepackage{bm}
\usepackage{color}
 \usepackage{multirow}
 \usepackage{rotating}

\input{epsf} 
\begin{document}

\title{
A Physics-Based Estimation of Mean Curvature Normal Vector for Triangulated Surfaces}
\author{Sudip Kumar Das$^1$, Mirza Cenanovic$^2$, and Junfeng Zhang$^1$\footnote{Corresponding author: Dr. Junfeng Zhang, Bharti School of Engineering, Laurentian University, 935 Ramsey Lake Road, Sudbury, ON P3E 2C6, Canada. Tel: 1-705-675-1151 ext. 2248; Fax: 1-705-675-4862; Email: jzhang@laurentian.ca.}
}
\affiliation{$^1$ Bharti School of Engineering, Laurentian University, 935 Ramsey Lake Road, Sudbury, Ontario, P3E 2C6, Canada; \\
$^2$ Department of Mechanical Engineering, Jonkoping University, SE-55111 Jonkoping, Sweden}
\date{\today}

\begin{abstract}
In this note, we derive an approximation for the mean curvature normal vector on vertices of triangulated surface meshes from the Young-Laplace equation and the force balance principle.
We then demonstrate that the approximation expression from our physics-based derivation is equivalent to the discrete Laplace-Beltrami operator approach in the literature.
This work, in addition to providing an alternative expression to calculate the mean curvature normal vector, can be further extended to other mesh structures, including non-triangular and heterogeneous meshes.
\end{abstract}


\maketitle

Estimating the normal and curvature for a triangulated surface is necessary for many computations and simulations, such as computer graphics, reverse engineering, medical image analysis, and multiphase flow simulations  
\cite{review1, review3, med_img, BIM_curvature, flow_review, WRJ_2013}.
Various methods have been proposed and readers can refer to several review articles
\cite{review1, review2, review3} and references therein.
Among them, Desbrun et al. \cite{Desbrun_1999, Meyer_Desbrun_2003} 
proposed an approximation method for differential attributes of triangular meshes.
Recognizing that the normal vector and mean curvature on a smooth surface can be 
expressed by the Laplace-Beltrami (LB) operator in differential geometry, 
the authors integrated the LB operator over a control surface around a vertex of the triangulated mesh for an estimation of the local mean curvature vector.
This method has since been employed in various applications, such as 
computer graphics \cite{CompGraph},
artificial intelligence \cite{AI},
biomedical engineering \cite{BioMedEng}.
In this note, we will present an alternative derivation for the discrete LB operator based on simple physics, namely, the Young–Laplace equation \cite{BIM_curvature, Neumann_book} and the force balance principle \cite{Statics_book}.
We will also show that our derived expression is equivalent to the discrete LB method by Desbrun et al. \cite{Desbrun_1999, Meyer_Desbrun_2003}.

\begin{figure}[t]
\begin{center}
\includegraphics[width= 1.0\textwidth]{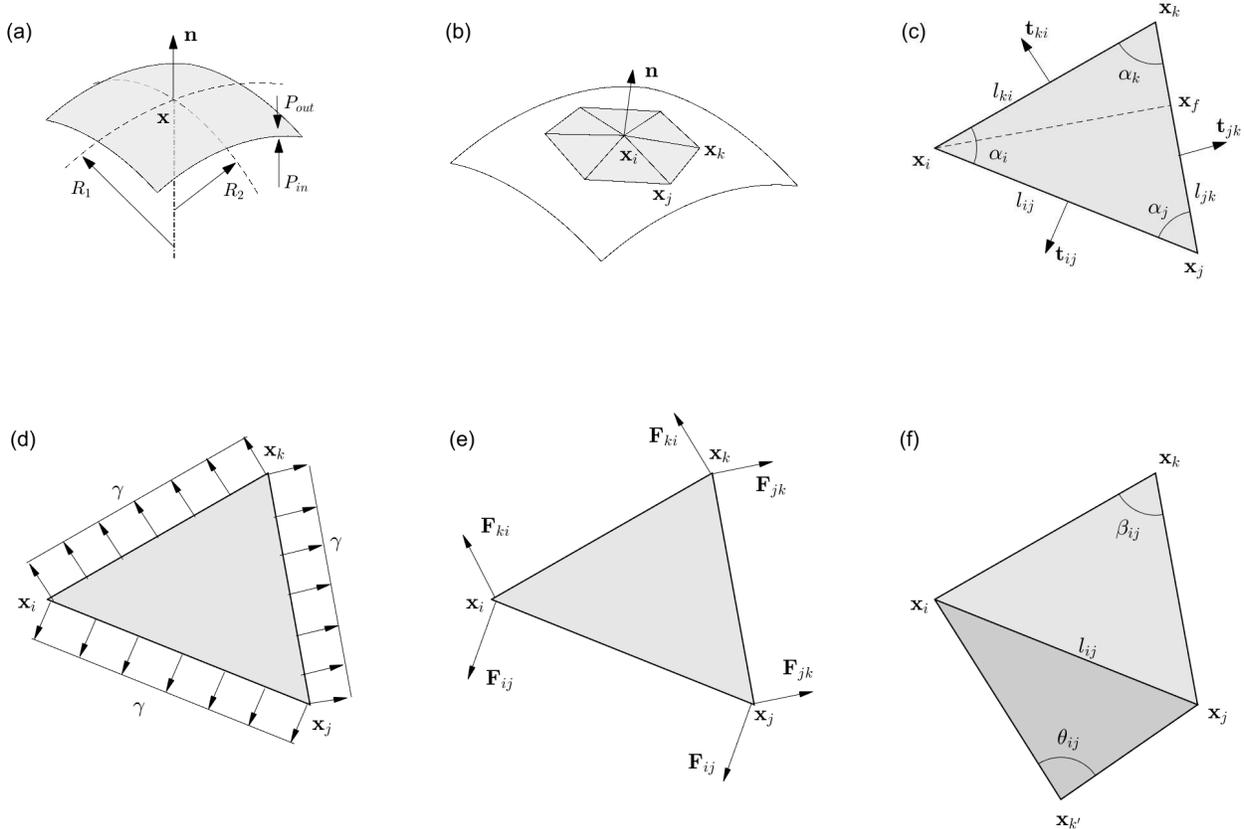}
\caption{Schematic drawings for the method description: (a) the original smooth surface; (b) a vertex ${\bf x}_i$ and its adjacent elements of the triangulated surface mesh; (c) notations for the interior angles, edge lengths, and outward edge perpendicular vectors of the element $E_{ijk}$; (d) the free-body diagram for the element $E_{ijk}$ with uniform force distributions along edges; (e) the replacement of the uniform force distributions with end forces; and (f) the two side elements $E_{ik^{\prime}j}$ and $E_{ijk}$ sharing the same edge and the two opposite angles $\beta_{ij}$ and $\theta_{ij}$.
}
\label{fig:schms}
\end{center}
\end{figure}

For a system consisting two static fluids separated by an interface, there may exist a difference in fluid pressure across the interface, and this pressure difference is given by the classical Young-Laplace equation \cite{Young_1805, Laplace_1805, BioFM, Neumann_book}:
\begin{equation}
\label{eq:YLE}
 P_{in}-P_{out}=\gamma \left( \frac{1}{R_1} + \frac{1}{R_2} \right)~,
\end{equation}
where $\gamma$ is the tension in the interface, $R_1$ and $R_2$ are the two principal radii of curvature,
and $P_{in}$ and $P_{out}$, respectively, are the fluid pressures inside and outside of the interface. 
The normal direction ${\bf  n}$ of the interface points into the outside fluid (Fig. 1a).
The mean curvature $H$ is the mean value of the principle curvatures $\kappa_1=1/R_1$ and $\kappa_2=1/R_2$, and thus Eq. (\ref{eq:YLE}) can be rewritten to
\begin{equation}
\label{eq:YLE2}
P_{in}-P_{out}=2 \gamma H~.
\end{equation}
The Young-Laplace equation Eq. (\ref{eq:YLE}) can be derived rigorously from physical principles such as the force balance and surface energy minimization \cite{Neumann_book, YLE_der}. 
On a differential surface area $\delta A$, over which the surface curvature and normal can be considered constant, the pressure difference $P_{in}-P_{out}$ generates a force
\begin{equation}
\label{eq:dF}
 \delta {\bf F} =(P_{in}-P_{out})  \delta A {\bf n} = 2\gamma H \delta A {\bf n}~.
\end{equation}

On the other hand, for a surface discretized into triangular elements (Fig. 1b), a nodal force 
${\bf F}_i$ is required on vertex ${\bf x}_i$ to maintain the structure configuration under the surface tension effect in triangular elements.
We consider one such adjacent element $E_{ijk}$ with vertices ${\bf x}_i$, ${\bf x}_j$ and ${\bf x}_k$,
and denote its interior angles as $\alpha$, edges as $l$, and the outward edge perpendicular vector  as ${\bf t}$, with corresponding subscripts (Fig. 1c).
A free-body diagram is constructed for this element Fig. 1d.
Here uniform force distributions of magnitude $\gamma$ in the outward perpendicular directions are applied along the element edges, according to the surface tension definition \cite{Neumann_book, Fluid_book}.
We then replace the uniform force distribution on each edge with two forces at the two end points of this edge \cite{Statics_book} (Fig. 1e):
\begin{equation}
F_{ij}=\frac{\gamma l_{ij}}{2}~,~~~
F_{jk}=\frac{\gamma l_{jk}}{2}~,~~~
F_{ki}=\frac{\gamma l_{kl}}{2}~,
\end{equation}
where $l_{ij}$, $l_{jk}$, and $l_{ki}$ are the edge lengths of the element with the subscripts indicating the two end nodes of the edge.
The nodal force ${\bf F}_{i, E_{ijk}}$, the force resulted at vertex ${\bf x}_i$ due to the surface tension $\gamma$ in element $E_{ijk}$, is
\begin{equation}
{\bf F}_{i,E_{ijk}}=F_{ki} {\bf t}_{ki}+ F_{ij} {\bf t}_{ij}~;
\end{equation}
and using some simple geometric relations one can  further simplify this expression to 
\begin{equation}
{\bf F}_{i,E_{ijk}}=-\gamma l_{jk} {\bf t}_{jk}~.
\label{eq:Fijk}
\end{equation}
It is interesting to see that the nodal force ${\bf F}_{i,E_{ijk}}$ at vortex ${\bf x}_i$ is simply the force required to counter balance the surface tension force $\gamma l_{jk}$ on the opposite edge.
The total force at vortex ${\bf x}_i$ is simply the sum of such nodal forces from all adjacent elements:
\begin{equation}
{\bf F}_i=-\gamma \sum_{(j,k)}l_{jk}{\bf t}_{jk}~.
\label{eq:Fi}
\end{equation}
Considering the triangulated mesh in Fig. 1b is a discrete approximation  of the original smooth surface in Fig. 1a, the nodal force ${\bf F}_i$ given in Eq. (\ref{eq:Fi}) is therefore an approximation of the force $\delta {\bf F}$ due to the pressure difference over a control volume surface patch $A_i$ associated with vortex ${\bf x}_i$
\begin{equation}
{\bf F}_i \approx  2\gamma H  A_i {\bf n}~.
\label{eq:FiP}
\end{equation}
Combining Eqs.  (\ref{eq:FiP}) and (\ref{eq:Fi}) and canceling the surface tension $\gamma$ yield the following approximation for the mean curvature normal $H{\bf n}$
\begin{equation}
\label{eq:Hn}
H{\bf n} \approx - \frac{1}{2A_i}\sum_{(j,k)}l_{jk}{\bf t}_{jk}~.
\end{equation}
With the nodal coordinates and connectivity available for the mesh network, the opposite edge length $l_{jk}$ and its outward direction ${\bf t}_{jk}$ in each adjacent element can be readily calculated. 
The control volume area $A_i$ can be calculated as the modified Voronoi area \cite{Meyer_Desbrun_2003}:
\begin{equation}
\label{eq:A-Meyer}
A_i=\sum_{(j,k)}A_{ijk}^v~~\mbox{with}~~~ A_{ijk}^v= 
\begin{cases}\frac{\cot \alpha_j |{\bf x}_k-{\bf x}_i|^2+
\cot \alpha_k |{\bf x}_j-{\bf x}_i|^2}{8}~, & \max(\alpha_i, \alpha_j, \alpha_k)\leq\frac{\pi}{2}~;\\
\frac{A_{ijk}}{2}~, & \alpha_i>\frac{\pi}{2}~;\\
\frac{A_ijk}{4}~,  & \alpha_j>\frac{\pi}{2}~\mbox{or}~
\alpha_k>\frac{\pi}{2}~,
\end{cases}
\end{equation}
where $A_{ijk}$ is the area of element $E_{ijk}$.

Next, we show that our approximation Eq. (\ref{eq:Hn}) actually is equivalent to that given by Desbrun et al. \cite{Desbrun_1999, Meyer_Desbrun_2003}.
Again we start with the element $E_{ijk}$ in Fig. 1c, and its inner angles are  $\alpha_i$, 
 $\alpha_j$, and $\alpha_k$ at vortice ${\bf x}_i$, ${\bf x}_j$ and ${\bf x}_k$, respectively.  
The position of the foot of altitude ${\bf x}_f$ on the base 
$|{\bf x}_k-{\bf x}_j|$ can be calculated from ${\bf x}_j$ and ${\bf x}_k$ with the distances as weight factors:
\begin{equation}
{\bf x}_f=\frac{|{\bf x}_k-{\bf x}_f|}{l_{jk}}{\bf x}_j +
\frac{|{\bf x}_f-{\bf x}_j|}{l_{jk}}{\bf x}_k~,
\end{equation}
and the outward vector ${\bf t}_{jk}$ is then
\begin{equation}
{\bf t}_{jk}=\frac{{\bf x}_f-{\bf x}_i}{|{\bf x}_f-{\bf x}_i|}
=\frac{
|{\bf x}_k-{\bf x}_f|{\bf x}_j +|{\bf x}_f-{\bf x}_j|{\bf x}_k
-l_{jk}{\bf x}_i}{l_{jk}|{\bf x}_f-{\bf x}_i|}
=\frac{\cot \alpha_k ({\bf x}_j-{\bf x}_i) +
\cot \alpha_j ({\bf x}_k-{\bf x}_i)} {l_{jk}}~,
\end{equation}
where the relationship $l_{jk}$=$|{\bf x}_f-{\bf x}_j|$ + $|{\bf x}_k-{\bf x}_f|$
has been utilized.
With this relationship obtained, we can rewrite Eq. (\ref{eq:Hn}) by
enumerating every link instead of every element connected at the same vertex ${\bf x}_i$ in the summation operation:
\begin{equation}
\label{eq:Hn2}
H{\bf n} ({\bf x}_i)\approx  \frac{1}{2A_i}\sum_j (\cot\beta_{ij}+\cot\theta_{ij})({\bf x}_i -{\bf x}_j)~.
\end{equation}
Here $\beta_{ij}$ and $\theta_{ij}$ are the opposite angles of the two elements sharing the edge $l_{ij}$ (Fig. 1f).
This is exactly the same expression given by Desbrun et al. \cite{Desbrun_1999, Meyer_Desbrun_2003}.
We do not present any numerical examples for our method Eq. (\ref{eq:Hn}) due to the mathematical identity to  Eq. (\ref{eq:Hn2}), which has been tested in numerous studies \cite{Desbrun_1999,Meyer_Desbrun_2003,  CompGraph, AI, BioMedEng}.

In summary, we have derived an estimation expression for the mean curvature vector at vertices of triangulated surfaces,
and shown that our expression is equivalent to the result from the discrete LB operator approach by Desbrun et al. \cite{Desbrun_1999, Meyer_Desbrun_2003}.
The key idea of our method can also be extended to other non-triangular and even heterogeneous surface meshes.
In such situations, instead of Eq. (\ref{eq:Fijk}), the finite element method \cite{FEM_book} can be utilized to calculate the nodal force due to the surface tension effect in each non-triangular surface element.

\section*{Acknowledgment}
This work was supported by 
the Natural Science and Engineering Research Council of Canada (NSERC). 
\bibliography{draft}

\begin{thebibliography}{20}
\expandafter\ifx\csname natexlab\endcsname\relax\def\natexlab#1{#1}\fi
\expandafter\ifx\csname bibnamefont\endcsname\relax
  \def\bibnamefont#1{#1}\fi
\expandafter\ifx\csname bibfnamefont\endcsname\relax
  \def\bibfnamefont#1{#1}\fi
\expandafter\ifx\csname citenamefont\endcsname\relax
  \def\citenamefont#1{#1}\fi
\expandafter\ifx\csname url\endcsname\relax
  \def\url#1{\texttt{#1}}\fi
\expandafter\ifx\csname urlprefix\endcsname\relax\def\urlprefix{URL }\fi
\providecommand{\bibinfo}[2]{#2}
\providecommand{\eprint}[2][]{\url{#2}}

\bibitem[{\citenamefont{Petitjean}(2002)}]{review1}
\bibinfo{author}{\bibfnamefont{S.}~\bibnamefont{Petitjean}},
  \bibinfo{journal}{{ACM} Computing Surveys} \textbf{\bibinfo{volume}{2}},
  \bibinfo{pages}{1} (\bibinfo{year}{2002}).

\bibitem[{\citenamefont{Nigam and Agrawal}(2013)}]{review3}
\bibinfo{author}{\bibfnamefont{S.}~\bibnamefont{Nigam}} \bibnamefont{and}
  \bibinfo{author}{\bibfnamefont{V.}~\bibnamefont{Agrawal}},
  \bibinfo{journal}{International Journal of Engineering Science and Innovative
  Technology} \textbf{\bibinfo{volume}{2}}, \bibinfo{pages}{330}
  (\bibinfo{year}{2013}).

\bibitem[{\citenamefont{Campbell and Summers}(2004)}]{med_img}
\bibinfo{author}{\bibfnamefont{S.~R.} \bibnamefont{Campbell}} \bibnamefont{and}
  \bibinfo{author}{\bibfnamefont{R.~M.} \bibnamefont{Summers}},
  \bibinfo{journal}{International Congress Series}
  \textbf{\bibinfo{volume}{1268}}, \bibinfo{pages}{999} (\bibinfo{year}{2004}).

\bibitem[{\citenamefont{Janssen and Anderson}(2007)}]{BIM_curvature}
\bibinfo{author}{\bibfnamefont{P.~J.~A.} \bibnamefont{Janssen}}
  \bibnamefont{and} \bibinfo{author}{\bibfnamefont{P.~D.}
  \bibnamefont{Anderson}}, \bibinfo{journal}{Physics of Fluids}
  \textbf{\bibinfo{volume}{19}}, \bibinfo{pages}{043602}
  (\bibinfo{year}{2007}).

\bibitem[{\citenamefont{Popinet}(2018)}]{flow_review}
\bibinfo{author}{\bibfnamefont{S.}~\bibnamefont{Popinet}},
  \bibinfo{journal}{Annual Review of Fluid Mechanics}
  \textbf{\bibinfo{volume}{50}}, \bibinfo{pages}{49} (\bibinfo{year}{2018}).

\bibitem[{\citenamefont{Wang}(2013)}]{WRJ_2013}
\bibinfo{author}{\bibfnamefont{R.}~\bibnamefont{Wang}},
  \bibinfo{journal}{Journal of Nanoparticle Research}
  \textbf{\bibinfo{volume}{15}}, \bibinfo{pages}{2128} (\bibinfo{year}{2013}).

\bibitem[{\citenamefont{Surazhsky et~al.}(2003)\citenamefont{Surazhsky, Magid,
  Soldea, Elber, and Rivlin}}]{review2}
\bibinfo{author}{\bibfnamefont{T.}~\bibnamefont{Surazhsky}},
  \bibinfo{author}{\bibfnamefont{E.}~\bibnamefont{Magid}},
  \bibinfo{author}{\bibfnamefont{O.}~\bibnamefont{Soldea}},
  \bibinfo{author}{\bibfnamefont{G.}~\bibnamefont{Elber}}, \bibnamefont{and}
  \bibinfo{author}{\bibfnamefont{E.}~\bibnamefont{Rivlin}}, in
  \emph{\bibinfo{booktitle}{Proceedings of the 2003 {IEEE} International
  Conference on Robotics and Automation}} (\bibinfo{address}{Taipei, Taiwan},
  \bibinfo{year}{2003}).

\bibitem[{\citenamefont{Desbrun et~al.}({1999})\citenamefont{Desbrun, Meyer,
  Schroder, and Barr}}]{Desbrun_1999}
\bibinfo{author}{\bibfnamefont{M.}~\bibnamefont{Desbrun}},
  \bibinfo{author}{\bibfnamefont{M.}~\bibnamefont{Meyer}},
  \bibinfo{author}{\bibfnamefont{P.}~\bibnamefont{Schroder}}, \bibnamefont{and}
  \bibinfo{author}{\bibfnamefont{A.~H.} \bibnamefont{Barr}}, in
  \emph{\bibinfo{booktitle}{Computer Graphics}} (\bibinfo{year}{{1999}}), pp.
  \bibinfo{pages}{{317--324}}.

\bibitem[{\citenamefont{Meyer et~al.}(2003)\citenamefont{Meyer, Desbrun, and
  Barr}}]{Meyer_Desbrun_2003}
\bibinfo{author}{\bibfnamefont{N.}~\bibnamefont{Meyer}},
  \bibinfo{author}{\bibfnamefont{M.}~\bibnamefont{Desbrun}}, \bibnamefont{and}
  \bibinfo{author}{\bibfnamefont{P.~S. A.~H.} \bibnamefont{Barr}}, in
  \emph{\bibinfo{booktitle}{Visualization and Mathematics}}, edited by
  \bibinfo{editor}{\bibfnamefont{H.~C.} \bibnamefont{Hege}} \bibnamefont{and}
  \bibinfo{editor}{\bibfnamefont{K.}~\bibnamefont{Polthier}}
  (\bibinfo{year}{2003}), vol.~\bibinfo{volume}{3}, pp.
  \bibinfo{pages}{35--57}.

\bibitem[{\citenamefont{Bickel et~al.}(2007)\citenamefont{Bickel, Botsch,
  Angst, Matusik, Otaduy, Pfister, and Gross}}]{CompGraph}
\bibinfo{author}{\bibfnamefont{B.}~\bibnamefont{Bickel}},
  \bibinfo{author}{\bibfnamefont{M.}~\bibnamefont{Botsch}},
  \bibinfo{author}{\bibfnamefont{R.}~\bibnamefont{Angst}},
  \bibinfo{author}{\bibfnamefont{W.}~\bibnamefont{Matusik}},
  \bibinfo{author}{\bibfnamefont{M.}~\bibnamefont{Otaduy}},
  \bibinfo{author}{\bibfnamefont{H.}~\bibnamefont{Pfister}}, \bibnamefont{and}
  \bibinfo{author}{\bibfnamefont{M.}~\bibnamefont{Gross}},
  \bibinfo{journal}{{ACM} Transactions on Graphics}
  \textbf{\bibinfo{volume}{26}}, \bibinfo{pages}{33} (\bibinfo{year}{2007}).

\bibitem[{\citenamefont{Lombaert et~al.}(2013)\citenamefont{Lombaert, Grady,
  Polimeni, and Cheriet}}]{AI}
\bibinfo{author}{\bibfnamefont{H.}~\bibnamefont{Lombaert}},
  \bibinfo{author}{\bibfnamefont{I.}~\bibnamefont{Grady}},
  \bibinfo{author}{\bibfnamefont{J.~R.} \bibnamefont{Polimeni}},
  \bibnamefont{and} \bibinfo{author}{\bibfnamefont{P.}~\bibnamefont{Cheriet}},
  \bibinfo{journal}{{IEEE} Transactions of Pattern Analysis and Machine
  Intelligence} \textbf{\bibinfo{volume}{35}}, \bibinfo{pages}{2143}
  (\bibinfo{year}{2013}).

\bibitem[{\citenamefont{Watson et~al.}(2017)\citenamefont{Watson, Sazonov,
  Zawieja, Moore, and {van Loon}}}]{BioMedEng}
\bibinfo{author}{\bibfnamefont{D.~J.} \bibnamefont{Watson}},
  \bibinfo{author}{\bibfnamefont{I.}~\bibnamefont{Sazonov}},
  \bibinfo{author}{\bibfnamefont{D.~C.} \bibnamefont{Zawieja}},
  \bibinfo{author}{\bibfnamefont{J.~E.} \bibnamefont{Moore}}, \bibnamefont{and}
  \bibinfo{author}{\bibfnamefont{R.}~\bibnamefont{{van Loon}}},
  \bibinfo{journal}{Journal of Biomechanics} \textbf{\bibinfo{volume}{64}},
  \bibinfo{pages}{172} (\bibinfo{year}{2017}).

\bibitem[{\citenamefont{Neumann et~al.}(2010)\citenamefont{Neumann, David, and
  Zou}}]{Neumann_book}
\bibinfo{author}{\bibfnamefont{A.~W.} \bibnamefont{Neumann}},
  \bibinfo{author}{\bibfnamefont{R.}~\bibnamefont{David}}, \bibnamefont{and}
  \bibinfo{author}{\bibfnamefont{Y.}~\bibnamefont{Zou}},
  \emph{\bibinfo{title}{Applied Surface Thermodynamics}}
  (\bibinfo{publisher}{{CRC} Press}, \bibinfo{address}{Boca Raton, US},
  \bibinfo{year}{2010}).

\bibitem[{\citenamefont{Hibbeler}(2016)}]{Statics_book}
\bibinfo{author}{\bibfnamefont{R.~C.} \bibnamefont{Hibbeler}},
  \emph{\bibinfo{title}{Engineering Mechanics: Statics}}
  (\bibinfo{publisher}{Pearson}, \bibinfo{year}{2016}).

\bibitem[{\citenamefont{Young}(1805)}]{Young_1805}
\bibinfo{author}{\bibfnamefont{T.}~\bibnamefont{Young}},
  \bibinfo{journal}{Philosophical Transactions of the Royal Society of London}
  \textbf{\bibinfo{volume}{95}}, \bibinfo{pages}{65} (\bibinfo{year}{1805}).

\bibitem[{\citenamefont{Laplace}(1805)}]{Laplace_1805}
\bibinfo{author}{\bibfnamefont{P.}~\bibnamefont{Laplace}},
  \emph{\bibinfo{title}{Supplément au dixième livre du Traité de Mécanique
  Céleste}} (\bibinfo{publisher}{Courcier}, \bibinfo{address}{Paris, France},
  \bibinfo{year}{1805}), vol.~\bibinfo{volume}{4}, p. \bibinfo{pages}{1–79}.

\bibitem[{\citenamefont{Rubenstein et~al.}(2012)\citenamefont{Rubenstein, Yin,
  and Frame}}]{BioFM}
\bibinfo{author}{\bibfnamefont{D.~A.} \bibnamefont{Rubenstein}},
  \bibinfo{author}{\bibfnamefont{W.}~\bibnamefont{Yin}}, \bibnamefont{and}
  \bibinfo{author}{\bibfnamefont{M.~D.} \bibnamefont{Frame}},
  \emph{\bibinfo{title}{Biofluid Mechanics}} (\bibinfo{publisher}{Elsevier},
  \bibinfo{address}{Oxford, UK}, \bibinfo{year}{2012}).

\bibitem[{\citenamefont{Siqveland and Skjæveland}(2014)}]{YLE_der}
\bibinfo{author}{\bibfnamefont{L.~M.} \bibnamefont{Siqveland}}
  \bibnamefont{and} \bibinfo{author}{\bibfnamefont{S.~M.}
  \bibnamefont{Skjæveland}} (\bibinfo{year}{2014}), \bibinfo{note}{dOI:
  10.13140/RG.2.1.4485.5768}.

\bibitem[{\citenamefont{Cengel and Cimbala}(2013)}]{Fluid_book}
\bibinfo{author}{\bibfnamefont{Y.~A.} \bibnamefont{Cengel}} \bibnamefont{and}
  \bibinfo{author}{\bibfnamefont{J.~M.} \bibnamefont{Cimbala}},
  \emph{\bibinfo{title}{Fluid Mechanics Fundamentals and Applications}}
  (\bibinfo{publisher}{McGraw-Hill}, \bibinfo{year}{2013}).

\bibitem[{\citenamefont{Whiteley}(2017)}]{FEM_book}
\bibinfo{author}{\bibfnamefont{J.}~\bibnamefont{Whiteley}},
  \emph{\bibinfo{title}{Finite Element Methods: A Practical Guide}}
  (\bibinfo{publisher}{Springer Nature}, \bibinfo{address}{Cham, Switzerland},
  \bibinfo{year}{2017}).

\end{thebibliography}
\end{document}